\documentclass[12pt]{myarticle}

\newtheorem{thm}{Theorem}

\newtheorem{definition}{Definition}
\newtheorem{ex}{Example}
\usepackage{tikz}
\usetikzlibrary{positioning}
\usepackage{subcaption}

\usepackage{txfonts}
\usepackage[colorlinks]{hyperref}

\newcommand{\n}{\noindent}
\usepackage{graphicx}
\begin{document}

\title{F-index of graphs based on new operations related to the join of graphs}

\author[ps]{Prosanta Sarkar}
\ead{prosantasarkar87@gmail.com}
\author[nd]{Nilanjan De}
\ead{de.nilanjan@rediffmail.com}
\author[ap]{Anita Pal}
\ead{anita.buie@gmail.com}

\address[ps;ap]{Department of Mathematics, National Institute of Technology, Durgapur, India.}

\address[nd]{Department of
Basic Sciences and Humanities (Mathematics),\\ Calcutta Institute of Engineering and Management, Kolkata, India.}

\begin{abstract}
The forgotten topological index or F-index of a graph is defined as the sum of degree cube of all the vertices of the graph. This index was introduced by Gutman and  Trinajesti\'{c} more than 40 years ago. In this paper, we derive F-index of new operations of subdivisions graphs related to the graphs join.
\end{abstract}
\begin{keyword}
F-index, $\mathcal{F}$-join of Graphs, Graph operations.
\end{keyword}

\maketitle

\section{Introduction}

Let $G$ be a simple graph with the vertex set $V(G)$ and edge set $E(G)$ respectively. The degree of a vertex $v$ in $G$ is the number of vertices in $G$ which are connected to $v$ by an edge and denoted by ${{d}_{G}}{(v)}$. A topological index is a graph invarient which is numerical parameter obtained from a graph which characterize its topology.Thus for two isomorphic graph $G$ and $H$ the value of a particular topological index must be same for both of them. In graph theory, there are many topological indices which have very useful applications in Chemistry, Bio-Chemistry, Molecular Biology, Nanotechnology for QSAR/QSPR investigation, isomer discrimination, pharmaceutical drug design and much more. The most thoroughly studied and oldest topological indices are first and second Zagreb indices which were introduced by Gutman and Trinajesti\'{c} \cite{gutm72} in a paper in 1972 to study structure dependency of the total $\pi$-electron energy$(\epsilon)$. They are respectively defined as
\[{{M}_{1}}(G)=\sum\limits_{v\in V(G)}{{{d}_{G}}{{(v)}^{2}}}=\sum\limits_{uv\in E(G)}{[{{d}_{G}}(u)+{{d}_{G}}(v)]}\] and \[{{M}_{2}}(G)=\sum\limits_{uv\in E(G)}{{{d}_{G}}(u){{d}_{G}}(v)}.\]
For further study we refer the reader to the following recent studies\cite{gutm04,khal09,mar14,das04,liu06,zho05}.

The ``forgotten topological index" was introduced in the same paper where first and second Zagreb index were introduced. But in 2015, Furtula and Gutman \cite{furt15} reinvestigated this index again and named this index as ``the forgotten topological index" or F-index. This index was defined as
\[{F}(G)=\sum\limits_{v\in V(G)}{{{d}_{G}}{{(v)}^{3}}}=\sum\limits_{uv\in E(G)}{[{{d}_{G}}(u)^2+{{d}_{G}}(v)^2]}.\]
We refer our reader to \cite{nde17,de15,nay16,de16,de18,den17} for some recent study and application of this index. The hyper Zagreb index, was introduced by Shirrdel et al. in \cite{shir13} and is defined as
\[{HM}(G)=\sum\limits_{uv\in E(G)}{[{{d}_{G}}(u)+{{d}_{G}}(v)]^2}.\]
For different mathematical and chemical studies of this index we refer our reader to \citep{ndeh17,gao17,vey16,far15,basav16,farah15,wan16}. One of the redefined version of the Zagreb index is defined as
\[{{ReZM}}(G)=\sum\limits_{uv\in E(G)}{{{d}_{G}}(u){{d}_{G}}(v)[{{d}_{G}}(u)+{{d}_{G}}(v)]}.\]
We refer the reader to \cite{ranj13,gao16,kum17,farah16} for some recent study of this redefined Zagreb index.

Note that, as usual, ${{P}_{n}}$ denotes a path with n vertices and $(n-1)$ edges whereas ${{C}_{n}}$ ($n \geqslant 3$) denotes a cycle graph with $n$ vertices. In this paper all over we deal with only simple and connected graphs.

\section{Preliminary}
Let  ${{G}_{1}}={(V({{G}_{1}}),E({{G}_{1}}))}$   and   ${{G}_{2}}={(V({{G}_{2}}),E({{G}_{2}}))}$   be two connected graphs such that $|V({{G}_{1}})|={{n}_{1}}$, $|V({{G}_{2}})|={{n}_{2}}$ and $|E({{G}_{1}})|={{m}_{1}}$, $|E({{G}_{2}})|={{m}_{2}}$ respectively. Different graph operations and derived graphs play an important role in graph theory. In this paper we consider one very important graph operations called join of graphs and some derived graphs such as different subdivision related graphs. Thus we first define them.

\begin{definition}
Joining of two graphs  ${{G}_{1}}$ and ${{G}_{2}}$,  denoted by  ${{G}_{1}}\vee{{G}_{2}}$, is an another graph which is obtained by joining each vertex of ${{G}_{1}}$ to all vertices of ${{G}_{2}}$,
such that \[{{d}_{{{G}_{1}}\vee{{G}_{2}}}}(v) = \left\{ \begin{array}{ll}
{{d}_{{{G}_{1}}}}(v)+{{n}_{2}},~if~v\in V({{G}_{1}})\\[2mm]
{{d}_{{{G}_{2}}}}(v)+{{n}_{1}},~if~v\in V({{G}_{2}}).\\[2mm]
\end{array}\right.\]
\end{definition}

\begin{definition}
S(G) is the graph which is obtained from $G$ by adding an extra vertex into each edge of $G$. In other words replaced each edge of $G$ by a path of length 2.
\end{definition}

\begin{definition}
The graph R(G) is obtained from $G$ by inserting an additional vertex into each edge of $G$ and joining each additional vertex to the end vertices of the corresponding edge of $G$.
\end{definition}

\begin{definition}
Q(G) is a graph derived from $G$ by adding a new vertex to each edge of $G$, then joining with edges those pairs of new vertices on adjacent edges of $G$.
\end{definition}

\begin{definition}
The total graph T(G) is derive from $G$ by inserting an new vertex to each edge of $G$, then joining each new vertex to the end vertices of the corresponding edge and joining with edges those pairs of new vertices on adjacent edges of $G$.
\end{definition}

Let $\mathcal{F}=\{S,R,Q,T\}$ and also let $I(G)$ denotes the set of vertices of $\mathcal{F}(G)$ which are inserted into each edge of G, so that $V{{({\mathcal{F}}{(G))}}}={{{V}{(G)}}\cup{{I}{(G)}}}$. Here we first define vertex $\mathcal{F}$-join and edge $\mathcal{F}$-join of graph operations based on the join of two connected graphs ${{G}_{1}}$ and ${{G}_{2}}$, which are defined as follows
\begin{definition}\cite{ndep17}
Let ${{G}_{1}}$ and ${{G}_{2}}$ be two simple graph the vertex $\mathcal{F}$-join graph of ${{G}_{1}}$ and ${{G}_{2}}$ is a graph obtained from $\mathcal{F}({{G}_{1}})$ and ${{G}_{2}}$ by joining each vertex of ${{G}_{1}}$ to all vertex of ${{G}_{2}}$, so that the vertex and edge sets are $V{(F({{G}_{1}}))}{\cup}V({{G}_{2}})$ and $E({{G}_{1}}){\cup}E({{G}_{2}}){\cup}[{{uv:u\in{V({{G}_{1}}),v\in{V({{G}_{2}})}}}}]$ respectively and is denoted by ${{G}_{1}}{\dot{\vee}_\mathcal{F}}{{G}_{2}}$.
\end{definition}
\begin{definition}\cite{ndep17}
Let ${{G}_{1}}$ and ${{G}_{2}}$ be two simple graph the edge $\mathcal{F}$-join graph of ${{G}_{1}}$ and ${{G}_{2}}$, is a graph obtained from $\mathcal{F}({{G}_{1}})$ and ${G}_{2}$ by joining each vertex of $I({{G}_{1}})$ to all vertex of ${{G}_{2}}$, so that the vertex and edge sets are $V{(\mathcal{F}({{G}_{1}}))}{\cup}V({{G}_{2}})$ and $E({{G}_{1}}){\cup}E({{G}_{2}}){\cup}[{{uv:u\in{I({{G}_{1}}),v\in{V({{G}_{2}})}}}}]$ respectively and is denoted by ${{G}_{1}}{\underline{\vee}_{\mathcal{F}}}{{G}_{2}}$.
\end{definition}

In this paper we will study the F-index of vertex
$\mathcal{F}$-join and edge $\mathcal{F}$-join of graphs related to different subdivision of graphs where $\mathcal{F}=\{S,R,Q,T\}$.

\section{Main Results:}
In the following sub section we consider F-index of vertex and edge $\mathcal{F}$-join of graphs for different values of $\mathcal{F}$ as $S,R,Q,T$ respectively. First we start with vertex and edge $S$-join of graph.

\subsection{Vertex and edge S-join of Graphs.}
In this section we derived closed formula of F-index for vertex and edge S-join of graphs and using these formula we obtained F-index of some special graphs. The figure of vertex and edge S-join of $P_3$ and $P_4$ are given in figure 1.

\begin{definition}
Let ${{G}_{1}}$ and ${{G}_{2}}$ be two vertex disjoint graphs the vertex S-join graph is obtained from $S({{G}_{1}})$ and ${G}_{2}$ by joining each vertex of $V({{G}_{1}})$ to every vertex of ${{G}_{2}}$ and is denoted by ${{G}_{1}}{\dot{\vee}_{S}}{{G}_{2}}$.
\end{definition}

\begin{thm}
If ${{G}_{1}}$ and ${{G}_{2}}$ be two connected graph then
\begin{eqnarray*}
{{F}}({{G}_{1}}{{\dot{\vee }}_{S}}{{G}_{2}})&=&{{F}}({{G}_{1}})+{{F}}({{G}_{2}})+3{{n}_{2}}{{M}_{1}}({{G}_{1}})+3{{n}_{1}}{{M}_{1}}({{G}_{2}})+6{{m}_{1}}{{n}_{2}^{2}}+6{{m}_{2}}{{n}_{1}^{2}}\\
          &&+n_{1}{{n}_{2}}({{n}_{1}}^2+n_{2}^{2})+8{{m}_{1}}.
\end{eqnarray*}
\end{thm}
\n\textit{Proof.} By definition of F-index, we have
\begin{eqnarray*}
{{F}}({{G}_{1}}{{\dot{\vee }}_{S}}{{G}_{2}})&=&\sum\limits_{v\in V({{G}_{1}}{{{\dot{\vee }}}_{S}}{{G}_{2}})}{{{d}_{{{G}_{1}}{{{\dot{\vee }}}_{S}}{{G}_{2}}}}{{(v)}^{3}}}\\
          &=&\sum\limits_{v\in V({{G}_{1}})}{{{d}_{{{G}_{1}}{{{\dot{\vee }}}_{S}}{{G}_{2}}}}{{(v)}^{3}}+\sum\limits_{v\in V({{G}_{2}})}{{{d}_{{{G}_{1}}{{{\dot{\vee }}}_{S}}{{G}_{2}}}}{{(v)}^{3}}+\sum\limits_{v\in I({{G}_{1}})}{{{d}_{{{G}_{1}}{{{\dot{\vee }}}_{S}}{{G}_{2}}}}{{(v)}^{3}}}}}\\
          &=&\sum\limits_{v\in V({{G}_{1}})}{{{({{d}_{{{G}_{1}}}}(v)+{{n}_{2}})}^{3}}+\sum\limits_{v\in V({{G}_{2}})}{{{({{d}_{{{G}_{2}}}}(v)+{{n}_{1}})}^{3}}+\sum\limits_{v\in I({{G}_{1}})}{{{2}^{3}}}}}\\
          &=&\sum\limits_{v\in V({{G}_{1}})}{[{{d}_{{{G}_{1}}}}{{(v)}^{3}}+3{{n}_{2}}{{d}_{{{G}_{1}}}}(v)^{2}+3{{n}_{2}}^{2}{{d}_{{{G}_{1}}}}(v)+n_{2}^{3}]}+\sum\limits_{v\in I({{G}_{1}})}{{{2}^{3}}}\\
          &&+\sum\limits_{v\in V({{G}_{2}})}{[{{d}_{{{G}_{2}}}}{{(v)}^{3}}+3{{n}_{1}}{{d}_{{{G}_{2}}}}(v)^{2}+3{{n}_{1}}^{2}{{d}_{{{G}_{2}}}}(v)+n_{1}^{3}]}\\
          &=&\sum\limits_{v\in V({{G}_{1}})}{{{d}_{{{G}_{1}}}}{{(v)}^{3}}}+3{{n}_{2}}\sum\limits_{v\in V({{G}_{1}})}{{{d}_{{{G}_{1}}}}(v)^{2}}+3{{n}_{2}}^{2}\sum\limits_{v\in V({{G}_{1}})}{{{d}_{{{G}_{1}}}}(v)}+{{n}_{1}}n_{2}^{3}\\
          &&+\sum\limits_{v\in V({{G}_{2}})}{{{d}_{{{G}_{2}}}}{{(v)}^{3}}}+3{{n}_{1}}\sum\limits_{v\in V({{G}_{2}})}{{{d}_{{{G}_{2}}}}(v)^{2}}+3{{n}_{1}}^{2}\sum\limits_{v\in V({{G}_{2}})}{{{d}_{{{G}_{2}}}}(v)}\\
          &&+{{n}_{2}}n_{1}^{3}+8{{m}_{1}}\\
          &=&{{F}}({{G}_{1}})+{{F}}({{G}_{2}})+3{{n}_{2}}{{M}_{1}}({{G}_{1}})+3{{n}_{1}}{{M}_{1}}({{G}_{2}})+6{{m}_{1}}{{n}_{2}^{2}}+6{{m}_{2}}{{n}_{1}^{2}}\\
          &&+n_{1}{{n}_{2}}({{n}_{1}}^2+n_{2}^{2})+8{{m}_{1}}.
\end{eqnarray*}
Which is the required result.      \qed
\begin{ex}
Using theorem 1, we get
\begin{eqnarray*}
(i)~{{F}}({{P}_{n}}{{\dot{\vee }}_{s}}{{P}_{m}})&=&(mn-6)(m^2+n^2)+6mn(m+n)+24mn-10m-2n-36,\\
(ii)~{{F}}({{P}_{n}}{{\dot{\vee }}_{s}}{{C}_{m}})&=&mn\{(m^2+n^2)+6(m+n)\}-6m^2+24mn-10m+16n-22,\\
(iii)~{{F}}({{C}_{n}}{{\dot{\vee }}_{s}}{{C}_{m}})&=&mn\{(m^2+n^2)+6(m+n)\}-6m^2+24mn+8m+16n,\\
(iv)~{{F}}({{C}_{n}}{{\dot{\vee }}_{s}}{{P}_{m}})&=&mn\{(m^2+n^2)+6(m+n)\}-6n^2+24mn+8m-2n-14.
\end{eqnarray*}
\end{ex}

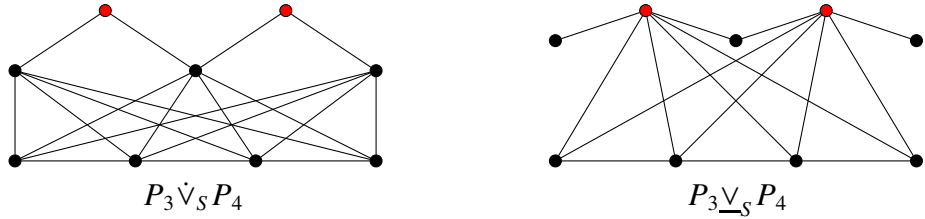
\begin{figure*}[h]
\begin{center}
\begin{tikzpicture}[scale=.4]
\tikzstyle{every node}=[draw, shape=circle, fill=red, scale=.4]
  \node (n1) at (6,7) {};
  \node (n2) at (12,7)  {};
\tikzstyle{every node}=[draw, shape=circle, fill=black, scale=.4]
  \node (n3) at (3,5)  {};
  \node (n4) at (9,5) {};
  \node (n5) at (15,5)  {};
  \node (n6) at (3,2)  {};
  \node (n7) at (7,2)  {};
  \node (n8) at (11,2)  {};
  \node (n9) at (15,2)   {};
  \foreach \from/\to in {n1/n3,n1/n4,n2/n4,n2/n5,n3/n6,n3/n7,n3/n8,n3/n9,n4/n6,n4/n7,n4/n8,n4/n9,n5/n6,n5/n7,n5/n8,n5/n9,n6/n7,n7/n8,n8/n9}
    \draw (\from) -- (\to);
\end{tikzpicture}
\hspace*{2cm}
\begin{tikzpicture}[scale=.4]
\tikzstyle{every node}=[draw, shape=circle, fill=red, scale=.4]
  \node (n1) at (6,7) {};
  \node (n2) at (12,7)  {};
\tikzstyle{every node}=[draw, shape=circle, fill=black, scale=.4]
  \node (n3) at (3,6)  {};
  \node (n4) at (9,6) {};
  \node (n5) at (15,6)  {};
  \node (n6) at (3,2)  {};
  \node (n7) at (7,2)  {};
  \node (n8) at (11,2)  {};
  \node (n9) at (15,2)   {};
  \foreach \from/\to in {n1/n3,n1/n4,n2/n4,n2/n5,n1/n6,n1/n7,n1/n8,n1/n9,n2/n6,n2/n7,n2/n8,n2/n9,n6/n7,n7/n8,n8/n9}
    \draw (\from) -- (\to);
\end{tikzpicture}
{${{P}_{3}}{\dot{\vee}_{S}}{{P}_{4}}$~~~~~~~~~~~~~~~~~~~~~~~~~~~~~~~~~~~~~~~~~~~~~~~~~~~~~~~~${{P}_{3}}{\underline{\vee}_{S}}{{P}_{4}}$}
\caption{The example of vertex S-join and edge S-join of graphs.}
\end{center}
\end{figure*}

\begin{definition}
Let ${{G}_{1}}$ and ${{G}_{2}}$ be two vertex disjoint graphs the edge S-join graph is obtained from $S({{G}_{1}})$ and ${G}_{2}$ by joining each vertex of $I({{G}_{1}})$ to every vertex of ${{G}_{2}}$ and is denoted by ${{G}_{1}}{\underline{\vee}_{S}}{{G}_{2}}$.
\end{definition}

\begin{thm}
If ${{G}_{1}}$ and ${{G}_{2}}$ be two connected graph then
\begin{eqnarray*}
{{F}}({{G}_{1}}{{\underline{\vee }}_{S}}{{G}_{2}})&=&{{F}}({{G}_{1}})+{{F}}({{G}_{2}})+3{{m}_{1}}{{M}_{1}}({{G}_{2}})+6{{m}_{1}^{2}}{{m}_{2}}+m_{1}({{n}_{2}+2)^3}+{{n}_{2}}{{{m}_{1}}}^{3}.
\end{eqnarray*}
\end{thm}
\n\textit{Proof.} From definition of F-index, we have
\begin{eqnarray*}
{{F}}({{G}_{1}}{{\underline{\vee }}_{S}}{{G}_{2}})&=&\sum\limits_{v\in V({{G}_{1}}{{{\underline{\vee }}}_{S}}{{G}_{2}})}{{{d}_{{{G}_{1}}{{{\underline{\vee }}}_{S}}{{G}_{2}}}}{{(v)}^{3}}}\\
          &=&\sum\limits_{v\in V({{G}_{1}})}{{{d}_{{{G}_{1}}{{{\underline{\vee }}}_{S}}{{G}_{2}}}}{{(v)}^{3}}+\sum\limits_{v\in V({{G}_{2}})}{{{d}_{{{G}_{1}}{{{\underline{\vee }}}_{S}}{{G}_{2}}}}{{(v)}^{3}}+\sum\limits_{v\in I({{G}_{1}})}{{{d}_{{{G}_{1}}{{{\underline{\vee }}}_{S}}{{G}_{2}}}}{{(v)}^{3}}}}}\\
          &=&\sum\limits_{v\in V({{G}_{1}})}{{{d}_{{{G}_{1}}}}(v)^{3}}+\sum\limits_{v\in V({{G}_{2}})}{{({{d}_{{{G}_{2}}}}(v)+{{m}_{1}})}^{3}}+\sum\limits_{v\in I({{G}_{1}})}{({{2}+{{n}_{2}}})^{3}}\\
          &=&\sum\limits_{v\in V({{G}_{1}})}{{{d}_{{{G}_{1}}}}{{(v)}^{3}}}+\sum\limits_{v\in I({{G}_{1}})}{({{2}+{{n}_{2}}})^{3}}\\
          &&+\sum\limits_{v\in V({{G}_{2}})}{[{{d}_{{{G}_{2}}}}{{(v)}^{3}}+3{{m}_{1}}{{d}_{{{G}_{2}}}}(v)^{2}+3{{m}_{1}}^{2}{{d}_{{{G}_{2}}}}(v)+m_{1}^{3}]}\\
          &=&\sum\limits_{v\in V({{G}_{1}})}{{{d}_{{{G}_{1}}}}{{(v)}^{3}}}+{{m}_{1}}{({{2}+{{n}_{2}}})^{3}}+\sum\limits_{v\in V({{G}_{2}})}{{{d}_{{{G}_{2}}}}{{(v)}^{3}}}\\
          &&+3{{m}_{1}}\sum\limits_{v\in V({{G}_{2}})}{{{d}_{{{G}_{2}}}}(v)^{2}}+3{{m}_{1}}^{2}\sum\limits_{v\in V({{G}_{2}})}{{{d}_{{{G}_{2}}}}(v)}+{{n}_{2}}m_{1}^{3}\\
          &=&{{F}}({{G}_{1}})+{{F}}({{G}_{2}})+3{{m}_{1}}{{M}_{1}}({{G}_{2}})+6{{m}_{1}^{2}}{{m}_{2}}\\
          &&+m_{1}({{n}_{2}+2)^3}+{{n}_{2}}{{{m}_{1}}}^{3}.
\end{eqnarray*}
Which is the desired result.                            \qed
\begin{ex}
Applying theorem 5, we get
\begin{eqnarray*}
(i)~{{F}}({{P}_{n}}{{\underline{\vee }}_{S}}{{P}_{m}})&=&(n-1)\{(m+2)^3+6(m-1)(n-1)+m(n-1)^2\}+12mn\\
&&-4m-10n-10,\\
(ii)~{{F}}({{P}_{n}}{{\underline{\vee }}_{S}}{{C}_{m}})&=&(n-1)\{(m+2)^3+6m(n-1)+m(n-1)^2\}+12mn\\
&&-4m+8n-14,\\
(iii)~{{F}}({{C}_{n}}{{\underline{\vee }}_{S}}{{C}_{m}})&=&n\{(m+2)^3+6mn+mn^2\}+12mn+8m+8n,\\
(iv)~{{F}}({{C}_{n}}{{\underline{\vee }}_{S}}{{P}_{m}})&=&n\{(m+2)^3+6n(m-1)+mn^2\}+12mn+8m-10n-14.
\end{eqnarray*}
\end{ex}

\subsection{Vertex and edge R-join of Graphs.}
Here we obtained F-index of vertex and edge R-join of graphs and applying these formula we derived F-index of some particular graphs. The figure of vertex and edge R-join of $P_3$ and $P_4$ are given in figure 2.

\begin{definition}
Let ${{G}_{1}}$ and ${{G}_{2}}$ be two vertex disjoint graphs the vertex R-join graph is obtained from $R({{G}_{1}})$ and ${G}_{2}$ by joining each vertex of $V({{G}_{1}})$ to every vertex of ${{G}_{2}}$ and is denoted by ${{G}_{1}}{\dot{\vee}_{R}}{{G}_{2}}$.
\end{definition}

\begin{thm}
If ${{G}_{1}}$ and ${{G}_{2}}$ be two connected graph then
\begin{eqnarray*}
{{F}}({{G}_{1}}{{\dot{\vee }}_{R}}{{G}_{2}})&=&8{{F}}({{G}_{1}})+{{F}}({{G}_{2}})+12{{n}_{2}}{{M}_{1}}({{G}_{1}})+3{{n}_{1}}{{M}_{1}}({{G}_{2}})+12{{m}_{1}}{{n}_{2}^{2}}+6{{m}_{2}}{{n}_{1}^{2}}\\
          &&+n_{1}{{n}_{2}}({{n}_{1}}^2+n_{2}^{2})+8{{m}_{1}}.
\end{eqnarray*}
\end{thm}
\n\textit{Proof.} We have from definition of F-index
\begin{eqnarray*}
{{F}}({{G}_{1}}{{\dot{\vee }}_{R}}{{G}_{2}})&=&\sum\limits_{v\in V({{G}_{1}}{{{\dot{\vee }}}_{R}}{{G}_{2}})}{{{d}_{{{G}_{1}}{{{\dot{\vee }}}_{R}}{{G}_{2}}}}{{(v)}^{3}}}\\
          &=&\sum\limits_{v\in V({{G}_{1}})}{{{d}_{{{G}_{1}}{{{\dot{\vee }}}_{R}}{{G}_{2}}}}{{(v)}^{3}}+\sum\limits_{v\in V({{G}_{2}})}{{{d}_{{{G}_{1}}{{{\dot{\vee }}}_{R}}{{G}_{2}}}}{{(v)}^{3}}+\sum\limits_{v\in I({{G}_{1}})}{{{d}_{{{G}_{1}}{{{\dot{\vee }}}_{R}}{{G}_{2}}}}{{(v)}^{3}}}}}\\
          &=&\sum\limits_{v\in V({{G}_{1}})}{{{(2{{d}_{{{G}_{1}}}}(v)+{{n}_{2}})}^{3}}+\sum\limits_{v\in V({{G}_{2}})}{{{({{d}_{{{G}_{2}}}}(v)+{{n}_{1}})}^{3}}+\sum\limits_{v\in I({{G}_{1}})}{{{2}^{3}}}}}\\
          &=&\sum\limits_{v\in V({{G}_{1}})}{[8{{d}_{{{G}_{1}}}}{{(v)}^{3}}+12{{n}_{2}}{{d}_{{{G}_{1}}}}(v)^{2}+6{{n}_{2}}^{2}{{d}_{{{G}_{1}}}}(v)+n_{2}^{3}]}+\sum\limits_{v\in I({{G}_{1}})}{{{2}^{3}}}\\
          &&+\sum\limits_{v\in V({{G}_{2}})}{[{{d}_{{{G}_{2}}}}{{(v)}^{3}}+3{{n}_{1}}{{d}_{{{G}_{2}}}}(v)^{2}+3{{n}_{1}}^{2}{{d}_{{{G}_{2}}}}(v)+n_{1}^{3}]}\\
          &=&8\sum\limits_{v\in V({{G}_{1}})}{{{d}_{{{G}_{1}}}}{{(v)}^{3}}}+12{{n}_{2}}\sum\limits_{v\in V({{G}_{1}})}{{{d}_{{{G}_{1}}}}(v)^{2}}+6{{n}_{2}}^{2}\sum\limits_{v\in V({{G}_{1}})}{{{d}_{{{G}_{1}}}}(v)}\\
          &&+{{n}_{1}}n_{2}^{3}+\sum\limits_{v\in V({{G}_{2}})}{{{d}_{{{G}_{2}}}}{{(v)}^{3}}}+3{{n}_{1}}\sum\limits_{v\in V({{G}_{2}})}{{{d}_{{{G}_{2}}}}(v)^{2}}\\
          &&+3{{n}_{1}}^{2}\sum\limits_{v\in V({{G}_{2}})}{{{d}_{{{G}_{2}}}}(v)}+{{n}_{2}}n_{1}^{3}+8{{m}_{1}}\\
          &=&8{{F}}({{G}_{1}})+{{F}}({{G}_{2}})+12{{n}_{2}}{{M}_{1}}({{G}_{1}})+3{{n}_{1}}{{M}_{1}}({{G}_{2}})+12{{m}_{1}}{{n}_{2}^{2}}+6{{m}_{2}}{{n}_{1}^{2}}\\
          &&+n_{1}{{n}_{2}}({{n}_{1}}^2+n_{2}^{2})+8{{m}_{1}}.
\end{eqnarray*}
which is the desired result.                                    \qed
\begin{ex}
From theorem 2, we get
\begin{eqnarray*}
(i)~{{F}}({{P}_{n}}{{\dot{\vee }}_{R}}{{P}_{m}})&=&mn(m^2+n^2+6n)+72mn-6n^2-76m\\
&&+54n-134,~ m,n\geq 2,\\
(ii)~{{F}}({{P}_{n}}{{\dot{\vee }}_{R}}{{C}_{m}})&=&mn(m^2+n^2+6n)+m^4+72mn-84m\\
&&+72n-120,~ m\geq 3,n\geq 2,\\
(iii)~{{F}}({{C}_{n}}{{\dot{\vee }}_{R}}{{C}_{m}})&=&mn(m^2+n^2+6n)+m^4+8n^4+72mn+8n,~ m,n\geq 3,\\
(iv)~{{F}}({{C}_{n}}{{\dot{\vee }}_{R}}{{P}_{m}})&=&mn(m^2+n^2)+6n^2(m-1)+72mn+16m\\
&&+46n-22,~ m\geq 3,n\geq 2.
\end{eqnarray*}
\end{ex}

\begin{figure*}[h]
\begin{center}
\begin{tikzpicture}[scale=.4]
\tikzstyle{every node}=[draw, shape=circle, fill=red, scale=.4]
  \node (n1) at (6,7) {};
  \node (n2) at (12,7)  {};
\tikzstyle{every node}=[draw, shape=circle, fill=black, scale=.4]
  \node (n3) at (3,5)  {};
  \node (n4) at (9,5) {};
  \node (n5) at (15,5)  {};
  \node (n6) at (3,2)  {};
  \node (n7) at (7,2)  {};
  \node (n8) at (11,2)  {};
  \node (n9) at (15,2)   {};
  \foreach \from/\to in {n1/n3,n1/n4,n2/n4,n2/n5,n3/n4,n4/n5,n3/n6,n3/n7,n3/n8,n3/n9,n4/n6,n4/n7,n4/n8,n4/n9,n5/n6,n5/n7,n5/n8,n5/n9,n6/n7,n7/n8,n8/n9}
    \draw (\from) -- (\to);
\end{tikzpicture}
\hspace*{2cm}
\begin{tikzpicture}[scale=.4]
\tikzstyle{every node}=[draw, shape=circle, fill=red, scale=.4]
  \node (n1) at (6,7) {};
  \node (n2) at (12,7)  {};
\tikzstyle{every node}=[draw, shape=circle, fill=black, scale=.4]
  \node (n3) at (3,6)  {};
  \node (n4) at (9,6) {};
  \node (n5) at (15,6)  {};
  \node (n6) at (3,2)  {};
  \node (n7) at (7,2)  {};
  \node (n8) at (11,2)  {};
  \node (n9) at (15,2)   {};
  \foreach \from/\to in {n1/n3,n1/n4,n2/n4,n2/n5,n3/n4,n4/n5,n1/n6,n1/n7,n1/n8,n1/n9,n2/n6,n2/n7,n2/n8,n2/n9,n6/n7,n7/n8,n8/n9}
    \draw (\from) -- (\to);
\end{tikzpicture}
{${{P}_{3}}{\dot{\vee}_{R}}{{P}_{4}}$~~~~~~~~~~~~~~~~~~~~~~~~~~~~~~~~~~~~~~~~~~~~~~~~~~~~~~~~${{P}_{3}}{\underline{\vee}_{R}}{{P}_{4}}$}
\caption{The example of vertex R-join and edge R-join of graphs.}
\end{center}
\end{figure*}
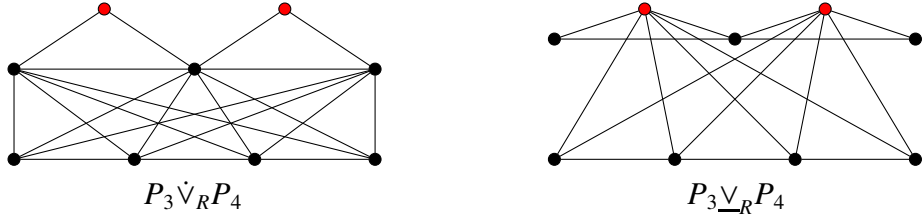

\begin{definition}
Let ${{G}_{1}}$ and ${{G}_{2}}$ be two vertex disjoint graphs the edge R-join graph is obtained from $R({{G}_{1}})$ and ${G}_{2}$ by joining each vertex of $I({{G}_{1}})$ to all vertex of ${{G}_{2}}$ and is denoted by ${{G}_{1}}{\underline{\vee}_{R}}{{G}_{2}}$.
\end{definition}

\begin{thm}
If ${{G}_{1}}$ and ${{G}_{2}}$ be two connected graph
\begin{eqnarray*}
{{F}}({{G}_{1}}{{\underline{\vee }}_{R}}{{G}_{2}})&=&8{{F}}({{G}_{1}})+{{F}}({{G}_{2}})+3{{m}_{1}}{{M}_{1}}({{G}_{2}})+6{{m}_{1}^{2}}{{m}_{2}}+m_{1}({{n}_{2}+2)^3}+{{n}_{2}}{{{m}_{1}}}^{3}.
\end{eqnarray*}
\end{thm}
\n\textit{Proof.} By definition of F-index, we have
\begin{eqnarray*}
{{F}}({{G}_{1}}{{\underline{\vee }}_{R}}{{G}_{2}})&=&\sum\limits_{v\in V({{G}_{1}}{{{\underline{\vee }}}_{R}}{{G}_{2}})}{{{d}_{{{G}_{1}}{{{\underline{\vee }}}_{R}}{{G}_{2}}}}{{(v)}^{3}}}\\
          &=&\sum\limits_{v\in V({{G}_{1}})}{{{d}_{{{G}_{1}}{{{\underline{\vee }}}_{R}}{{G}_{2}}}}{{(v)}^{3}}+\sum\limits_{v\in V({{G}_{2}})}{{{d}_{{{G}_{1}}{{{\underline{\vee }}}_{R}}{{G}_{2}}}}{{(v)}^{3}}+\sum\limits_{v\in I({{G}_{1}})}{{{d}_{{{G}_{1}}{{{\underline{\vee }}}_{R}}{{G}_{2}}}}{{(v)}^{3}}}}}\\
          &=&\sum\limits_{v\in V({{G}_{1}})}{{(2{{d}_{{{G}_{1}}}}(v))}^{3}}+\sum\limits_{v\in V({{G}_{2}})}{{({{d}_{{{G}_{2}}}}(v)+{{m}_{1}})}^{3}}+\sum\limits_{v\in I({{G}_{1}})}{({{2}+{{n}_{2}}})^{3}}\\
          &=&8\sum\limits_{v\in V({{G}_{1}})}{{{d}_{{{G}_{1}}}}{{(v)}^{3}}}+\sum\limits_{v\in I({{G}_{1}})}{({{2}+{{n}_{2}}})^{3}}\\
          &&+\sum\limits_{v\in V({{G}_{2}})}{[{{d}_{{{G}_{2}}}}{{(v)}^{3}}+3{{m}_{1}}{{d}_{{{G}_{2}}}}(v)^{2}+3{{m}_{1}}^{2}{{d}_{{{G}_{2}}}}(v)+m_{1}^{3}]}\\
          &=&8\sum\limits_{v\in V({{G}_{1}})}{{{d}_{{{G}_{1}}}}{{(v)}^{3}}}+{{m}_{1}}{({{2}+{{n}_{2}}})^{3}}+\sum\limits_{v\in V({{G}_{2}})}{{{d}_{{{G}_{2}}}}{{(v)}^{3}}}\\
          &&+3{{m}_{1}}\sum\limits_{v\in V({{G}_{2}})}{{{d}_{{{G}_{2}}}}(v)^{2}}+3{{m}_{1}}^{2}\sum\limits_{v\in V({{G}_{2}})}{{{d}_{{{G}_{2}}}}(v)}+{{n}_{2}}m_{1}^{3}\\
          &=&8{{F}}({{G}_{1}})+{{F}}({{G}_{2}})+3{{m}_{1}}{{M}_{1}}({{G}_{2}})+6{{m}_{1}^{2}}{{m}_{2}}\\
          &&+m_{1}({{n}_{2}+2)^3}+{{n}_{2}}{{{m}_{1}}}^{3}.
\end{eqnarray*}
Which is the required result.                         \qed
\begin{ex}
Applying theorem 6, we get
\begin{eqnarray*}
(i)~{{F}}({{P}_{n}}{{\underline{\vee }}_{R}}{{P}_{m}})&=&(n-1)(m+2)^3+m(n-1)^3+6(m-1)(n-1)^2+12mn\\
&&-4m+46n-94,~m,n\geq 2,\\
(ii)~{{F}}({{P}_{n}}{{\underline{\vee }}_{R}}{{C}_{m}})&=&(n-1)(m+2)^3+m(n-1)^2(n+5)+12mn-4m\\
&&+64n-112,~m\geq 3,n\geq 2,\\
(iii)~{{F}}({{C}_{n}}{{\underline{\vee }}_{R}}{{C}_{m}})&=&n(m+2)^3+mn^3+6mn^2+12mn+8m+64n,~m,n\geq 3,\\
(iv)~{{F}}({{C}_{n}}{{\underline{\vee }}_{R}}{{P}_{m}})&=&n(m+2)^3+mn^3+6(m-1)n^2+12mn+8m\\
&&+46n-14,~m\geq 2,n\geq 3.
\end{eqnarray*}
\end{ex}
\subsection{Vertex and edge Q-join of Graphs.}
In this section we obtained the F-index of vertex and edge Q-join of graphs and using these results
we obtained F-index of some special classes of graphs. The figure of vertex and edge Q-join of $P_3$ and $P_4$ are given in figure 3. Let us denote $M_4(G)$ as follows
\[{{M}_4}(G)=\sum\limits_{v\in V(G)}{{{d}_{G}}{{(v)}^4}}=\sum\limits_{uv\in E(G)}{[{{d}_{G}}(u)^3+{{d}_{G}}(v)^3]}\]
which is a particular case of general Zagreb index. In the following we use this topological index to represents vertex and edge Q-join of graphs.

\begin{definition}
Let ${{G}_{1}}$ and ${{G}_{2}}$ be two vertex disjoint graphs the vertex Q-join graph is obtained from $Q({{G}_{1}})$ and ${G}_{2}$ by joining each vertex of $V({{G}_{1}})$ to every vertex of ${{G}_{2}}$ and is denoted by ${{G}_{1}}{\dot{\vee}_{Q}}{{G}_{2}}$.
\end{definition}

\begin{thm}
If ${{G}_{1}}$ and ${{G}_{2}}$ be two connected graph then
\begin{eqnarray*}
{{F}}(G{}_{1}{{\dot{\vee }}_{Q}}G{}_{2})&=&{{F}}({{G}_{1}})+{{F}}({{G}_{2}})+3{{n}_{2}}{{M}_{1}}({{G}_{1}})+3{{n}_{1}}{{M}_{1}}({{G}_{2}})+{{M}_{4}}({{G}_{1}})\\
          &&+3{{ReZM}({{G}_{1}}})+6{{m}_{1}}{{n}_{2}^{2}}+6{{m}_{2}}{{n}_{1}^{2}}+n_{1}{{n}_{2}}({{n}_{1}}^2+n_{2}^{2}).
\end{eqnarray*}
\end{thm}
\n\textit{Proof.} From definition of F-index, we have
\begin{eqnarray*}
{{F}}({{G}_{1}}{{\dot{\vee }}_{Q}}{{G}_{2}})&=&\sum\limits_{v\in V({{G}_{1}}{{{\dot{\vee }}}_{Q}}{{G}_{2}})}{{{d}_{{{G}_{1}}{{{\dot{\vee }}}_{Q}}{{G}_{2}}}}{{(v)}^{3}}}\\
          &=&\sum\limits_{v\in V({{G}_{1}})}{{{d}_{{{G}_{1}}{{{\dot{\vee }}}_{Q}}{{G}_{2}}}}{{(v)}^{3}}+\sum\limits_{v\in V({{G}_{2}})}{{{d}_{{{G}_{1}}{{{\dot{\vee }}}_{Q}}{{G}_{2}}}}{{(v)}^{3}}+\sum\limits_{v\in I({{G}_{1}})}{{{d}_{{{G}_{1}}{{{\dot{\vee }}}_{Q}}{{G}_{2}}}}{{(v)}^{3}}}}}\\
          &=&\sum\limits_{v\in V({{G}_{1}})}{{{({{d}_{{{G}_{1}}}}(v)+{{n}_{2}})}^{3}}}+\sum\limits_{v\in V({{G}_{2}})}{{{({{d}_{{{G}_{2}}}}(v)+{{n}_{1}})}^{3}}}\\
          &&+\sum\limits_{uv\in E({{G}_{1}})}({{{{d}_{{{G}_{1}}}}(u)+{{d}_{{{G}_{1}}}}(v)})^{3}}\\
          &=&\sum\limits_{v\in V({{G}_{1}})}{[{{d}_{{{G}_{1}}}}{{(v)}^{3}}+3{{n}_{2}}{{d}_{{{G}_{1}}}}(v)^{2}+3{{n}_{2}}^{2}{{d}_{{{G}_{1}}}}(v)+n_{2}^{3}]}\\
          &&+\sum\limits_{v\in V({{G}_{2}})}{[{{d}_{{{G}_{2}}}}{{(v)}^{3}}+3{{n}_{1}}{{d}_{{{G}_{2}}}}(v)^{2}+3{{n}_{1}}^{2}{{d}_{{{G}_{2}}}}(v)+n_{1}^{3}]}\\
          &&+\sum\limits_{uv\in E({{G}_{1}})}[{{d}_{{{G}_{1}}}}(u)^{3}+{{d}_{{{G}_{1}}}}(v)^{3}+3{{d}_{{{G}_{1}}}}(u){{d}_{{{G}_{1}}}}(v)({{d}_{{{G}_{1}}}}(u)+{{d}_{{{G}_{1}}}}(v))]\\
          &=&\sum\limits_{v\in V({{G}_{1}})}{{{d}_{{{G}_{1}}}}{{(v)}^{3}}}+3{{n}_{2}}\sum\limits_{v\in V({{G}_{1}})}{{{d}_{{{G}_{1}}}}(v)^{2}}+3{{n}_{2}}^{2}\sum\limits_{v\in V({{G}_{1}})}{{{d}_{{{G}_{1}}}}(v)}+{{n}_{1}}n_{2}^{3}\\
          &&+\sum\limits_{v\in V({{G}_{2}})}{{{d}_{{{G}_{2}}}}{{(v)}^{3}}}+3{{n}_{1}}\sum\limits_{v\in V({{G}_{2}})}{{{d}_{{{G}_{2}}}}(v)^{2}}+3{{n}_{1}}^{2}\sum\limits_{v\in V({{G}_{2}})}{{{d}_{{{G}_{2}}}}(v)}\\
          &&+{{n}_{2}}n_{1}^{3}+\sum\limits_{uv\in E({{G}_{1}})}({{d}_{{{G}_{1}}}}(u)^{3}+{{d}_{{{G}_{1}}}}(v)^{3})\\
          &&+3\sum\limits_{uv\in E({{G}_{1}})}{{d}_{{{G}_{1}}}}(u){{d}_{{{G}_{1}}}}(v)({{d}_{{{G}_{1}}}}(u)+{{d}_{{{G}_{1}}}}(v))\\
          &=&{{F}}({{G}_{1}})+{{F}}({{G}_{2}})+3{{n}_{2}}{{M}_{1}}({{G}_{1}})+3{{n}_{1}}{{M}_{1}}({{G}_{2}})+{{M}_{4}}({{G}_{1}})\\
          &&+3{{ReZM}({{G}_{1}}})+6{{m}_{1}}{{n}_{2}^{2}}+6{{m}_{2}}{{n}_{1}^{2}}+n_{1}{{n}_{2}}({{n}_{1}}^2+n_{2}^{2}).
\end{eqnarray*}
Which is the required result.                                        \qed
\begin{ex}
Using theorem 3, we get
\begin{eqnarray*}
(i)~{{F}}({{P}_{n}}{{\dot{\vee }}_{Q}}{{P}_{m}})&=&mn(m^2+n^2)+6m^2(n-1)+6n^2(m-1)+24mn\\
&&-10m+54n-166,~ m,n\geq 3,\\
(ii)~{{F}({{P}_{n}}{{\dot{\vee }}_{Q}}{{C}_{m}})}&=&mn(m^2+n^2)+6m^2(n-1)+6n^2m+24mn\\
&&-10m+72n-152,~ m,n\geq 3 ,\\
(iii)~{{F}}({{C}_{n}}{{\dot{\vee }}_{Q}}{{C}_{m}})&=&mn(m^2+n^2)+6mn(m+n)+24mn+8m+72n,  ~  m,n\geq 3,\\
(iv)~{{F}({{C}_{n}}{{\dot{\vee }}_{Q}}{{P}_{m}})}&=&mn(m^2+n^2)+6mn(m+n)-6n^2+24mn+8m\\
&&+54n-14, ~ m,n\geq 3.
\end{eqnarray*}
\end{ex}

\begin{figure*}[h]
\begin{center}
\begin{tikzpicture}[scale=.4]
\tikzstyle{every node}=[draw, shape=circle, fill=red, scale=.4]
  \node (n1) at (6,7) {};
  \node (n2) at (12,7)  {};
\tikzstyle{every node}=[draw, shape=circle, fill=black, scale=.4]
  \node (n3) at (3,5)  {};
  \node (n4) at (9,5) {};
  \node (n5) at (15,5)  {};
  \node (n6) at (3,2)  {};
  \node (n7) at (7,2)  {};
  \node (n8) at (11,2)  {};
  \node (n9) at (15,2)   {};
  \foreach \from/\to in {n1/n3,n1/n4,n2/n4,n2/n5,n1/n2,n3/n6,n3/n7,n3/n8,n3/n9,n4/n6,n4/n7,n4/n8,n4/n9,n5/n6,n5/n7,n5/n8,n5/n9,n6/n7,n7/n8,n8/n9}
    \draw (\from) -- (\to);
\end{tikzpicture}
\hspace*{2cm}
\begin{tikzpicture}[scale=.4]
\tikzstyle{every node}=[draw, shape=circle, fill=red, scale=.4]
  \node (n1) at (6,7) {};
  \node (n2) at (12,7)  {};
\tikzstyle{every node}=[draw, shape=circle, fill=black, scale=.4]
  \node (n3) at (3,6)  {};
  \node (n4) at (9,6) {};
  \node (n5) at (15,6)  {};
  \node (n6) at (3,2)  {};
  \node (n7) at (7,2)  {};
  \node (n8) at (11,2)  {};
  \node (n9) at (15,2)   {};
  \foreach \from/\to in {n1/n3,n1/n4,n2/n4,n2/n5,n1/n2,n1/n6,n1/n7,n1/n8,n1/n9,n2/n6,n2/n7,n2/n8,n2/n9,n6/n7,n7/n8,n8/n9}
    \draw (\from) -- (\to);
\end{tikzpicture}
{${{P}_{3}}{\dot{\vee}_{Q}}{{P}_{4}}$~~~~~~~~~~~~~~~~~~~~~~~~~~~~~~~~~~~~~~~~~~~~~~~~~~~~~~~${{P}_{3}}{\underline{\vee}_{Q}}{{P}_{4}}$}
\caption{The example of vertex Q-join and edge Q-join of graphs.}
\end{center}
\end{figure*}
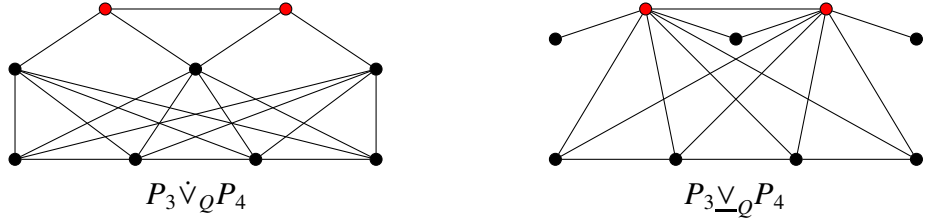

\begin{definition}
Let ${{G}_{1}}$ and ${{G}_{2}}$ be two vertex disjoint graphs the edge Q-join graph is obtained from $Q({{G}_{1}})$ and ${G}_{2}$ by joining each vertex of $I({{G}_{1}})$ to every vertex of ${{G}_{2}}$ and is denoted by ${{G}_{1}}{\underline{\vee}_{Q}}{{G}_{2}}$.
\end{definition}

\begin{thm}
If ${{G}_{1}}$ and ${{G}_{2}}$ be two connected graph then
\begin{eqnarray*}
{F}({{G}_{1}}{{\underline{\vee }}_{Q}}{{G}_{2}})&=&{{F}}({{G}_{1}})+{{F}}({{G}_{2}})+3{{n}_{2}}^{2}{{M}_{1}}({{G}_{1}})+3{{m}_{1}}{{M}_{1}}({{G}_{2}})+{{M}_{4}}({{G}_{1}})\\
          &&+3{{n}_{2}}{{HM}({{G}_{1}})}+3{{ReZM}({{G}_{1}})}+m_{1}^{2}(6{{m}_{2}}+{{m}_{1}}{{n}_{2}})+{{m}_{1}}{{{n}_{2}}}^{3}.
\end{eqnarray*}
\end{thm}
\n\textit{Proof.} We have from definition of F-index
\begin{eqnarray*}
{{F}}({{G}_{1}}{{\underline{\vee }}_{Q}}{{G}_{2}})&=&\sum\limits_{v\in V({{G}_{1}}{{{\underline{\vee }}}_{Q}}{{G}_{2}})}{{{d}_{{{G}_{1}}{{{\underline{\vee }}}_{Q}}{{G}_{2}}}}{{(v)}^{3}}}\\
          &=&\sum\limits_{v\in V({{G}_{1}})}{{{d}_{{{G}_{1}}{{{\underline{\vee }}}_{Q}}{{G}_{2}}}}{{(v)}^{3}}+\sum\limits_{v\in V({{G}_{2}})}{{{d}_{{{G}_{1}}{{{\underline{\vee }}}_{Q}}{{G}_{2}}}}{{(v)}^{3}}+\sum\limits_{v\in I({{G}_{1}})}{{{d}_{{{G}_{1}}{{{\underline{\vee }}}_{Q}}{{G}_{2}}}}{{(v)}^{3}}}}}\\
          &=&\sum\limits_{v\in V({{G}_{1}})}{{{d}_{{{G}_{1}}}}(v)^{3}}+\sum\limits_{v\in V({{G}_{2}})}{{({{d}_{{{G}_{2}}}}(v)+{{m}_{1}})}^{3}}\\
          &&+\sum\limits_{uv\in E({{G}_{1}})}({{{{d}_{{{G}_{1}}}}(u)+{{d}_{{{G}_{1}}}}(v)+{{n}_{2}})^{3}}}\\
          &=&\sum\limits_{v\in V({{G}_{1}})}{{{d}_{{{G}_{1}}}}{{(v)}^{3}}}+\sum\limits_{v\in V({{G}_{2}})}{[{{d}_{{{G}_{2}}}}{{(v)}^{3}}+3{{m}_{1}}{{d}_{{{G}_{2}}}}(v)^{2}+3{{m}_{1}}^{2}{{d}_{{{G}_{2}}}}(v)}\\
          &&+m_{1}^{3}]+\sum\limits_{uv\in E({{G}_{1}})}[({{{{d}_{{{G}_{1}}}}(u)+{{d}_{{{G}_{1}}}}(v))^{3}}}+3{{n}_{2}}{{({{d}_{{{G}_{1}}}}(u)+{{d}_{{{G}_{1}}}}(v))^{2}}}\\
          &&+3{{n}_{2}}^{2}{{({{d}_{{{G}_{1}}}}(u)+{{d}_{{{G}_{1}}}}(v))}}+{{n}_{2}}^{3}]\\
          &=&\sum\limits_{v\in V({{G}_{1}})}{{{d}_{{{G}_{1}}}}{{(v)}^{3}}}+\sum\limits_{v\in V({{G}_{2}})}{{{d}_{{{G}_{2}}}}{{(v)}^{3}}}+3{{m}_{1}}\sum\limits_{v\in V({{G}_{2}})}{{{d}_{{{G}_{2}}}}(v)^{2}}\\
          &&+3{{m}_{1}}^{2}\sum\limits_{v\in V({{G}_{2}})}{{{d}_{{{G}_{2}}}}(v)}+{{n}_{2}}m_{1}^{3}+\sum\limits_{uv\in E({{G}_{1}})}({{d}_{{{G}_{1}}}}(u)^{3}+{{d}_{{{G}_{1}}}}(v)^{3})\\
          &&+3\sum\limits_{uv\in E({{G}_{1}})}{{d}_{{{G}_{1}}}}(u){{d}_{{{G}_{1}}}}(v)({{d}_{{{G}_{1}}}}(u)+{{d}_{{{G}_{1}}}}(v))+{{n}_{2}}^{3}m_{1}\\
          &&+3{{n}_{2}}\sum\limits_{uv\in E({{G}_{1}})}({{d}_{{{G}_{1}}}}(u)+{{d}_{{{G}_{1}}}}(v))^{2}+3{{n}_{2}}^{2}\sum\limits_{uv\in E({{G}_{1}})}({{d}_{{{G}_{1}}}}(u)+{{d}_{{{G}_{1}}}}(v))\\
          &=&{{F}}({{G}_{1}})+{{F}}({{G}_{2}})+3{{n}_{2}}^{2}{{M}_{1}}({{G}_{1}})+3{{m}_{1}}{{M}_{1}}({{G}_{2}})+{{M}_{4}}({{G}_{1}})\\
          &&+3{{n}_{2}}{{HM}({{G}_{1}})}+3{{ReZM}({{G}_{1}})}+m_{1}^{2}(6{{m}_{2}}+{{m}_{1}}{{n}_{2}})+{{m}_{1}}{{{n}_{2}}}^{3}.
\end{eqnarray*}
Which is the desired result.                                      \qed
\begin{ex}
Using theorem 7, we get
\begin{eqnarray*}
(i)~{{F}}({{P}_{n}}{{\underline{\vee }}_{Q}}{{P}_{m}})&=&m(n-1)\{(n-1)^2+m^2\}+3m^2(4n-6)+6(m-1)(n-1)^2\\
&&+60mn-94m+54n-148,~m\geq 3,n\ge 4,\\
(ii)~{{F}}({{P}_{n}}{{\underline{\vee }}_{Q}}{{C}_{m}})&=&m(n-1)\{(n-1)^2+m^2\}+3m^2(4n-6)+6m(n-1)^2\\
&&+60mn-94m+72n-152,~m\geq 3,n\ge 4,\\
(iii)~{{F}}({{C}_{n}}{{\underline{\vee }}_{Q}}{{C}_{m}})&=&mn(m^2+n^2)+12m^2n+6mn^2+60mn+8m\\
&&+72n,~m\geq 3,n\ge 4,\\
(iv)~{{F}}({{C}_{n}}{{\underline{\vee }}_{Q}}{{P}_{m}})&=&mn(m^2+n^2)+12m^2n+6mn^2-6n^2+60mn+8m\\
&&+6n-14,~m\geq 3,n\ge 4.
\end{eqnarray*}
\end{ex}

\subsection{Vertex and edge T-join of Graphs.}
Here we derived the closed formula of F-index for vertex and edge T-join of graphs and applying these formula we obtained the F-index of some particular graphs. The figure of vertex and edge T-join of $P_3$ and $P_4$ are given in figure 4.
\begin{definition}
Let ${{G}_{1}}$ and ${{G}_{2}}$ be two vertex disjoint graphs the vertex T-join graph is obtained from $T({{G}_{1}})$ and ${G}_{2}$ by joining each vertex of $V({{G}_{1}})$ to every vertex of ${{G}_{2}}$ and is denoted by ${{G}_{1}}{\dot{\vee}_{T}}{{G}_{2}}$.
\end{definition}

\begin{thm}
If ${{G}_{1}}$ and ${{G}_{2}}$ be two connected graph then
\begin{eqnarray*}
{{F}}({{G}_{1}}{{\dot{\vee }}_{T}}{{G}_{2}})&=&8{{F}}({{G}_{1}})+{{F}}({{G}_{2}})+12{{n}_{2}}{{M}_{1}}({{G}_{1}})+3{{n}_{1}}{{M}_{1}}({{G}_{2}})+{{M}_{4}}({{G}_{1}})\\
          &&+3{{ReZM}({{G}_{1}}})+12{{m}_{1}}{{n}_{2}^{2}}+6{{m}_{2}}{{n}_{1}^{2}}+n_{1}{{n}_{2}}({{n}_{1}}^2+n_{2}^{2}).
\end{eqnarray*}
\end{thm}
\n\textit{Proof.} By definition of F-index, we have
\begin{eqnarray*} 	
{{F}}({{G}_{1}}{{\dot{\vee }}_{T}}{{G}_{2}})&=&\sum\limits_{v\in V({{G}_{1}}{{{\dot{\vee }}}_{T}}{{G}_{2}})}{{{d}_{{{G}_{1}}{{{\dot{\vee }}}_{T}}{{G}_{2}}}}{{(v)}^{3}}}\\
          &=&\sum\limits_{v\in V({{G}_{1}})}{{{d}_{{{G}_{1}}{{{\dot{\vee }}}_{T}}{{G}_{2}}}}{{(v)}^{3}}+\sum\limits_{v\in V({{G}_{2}})}{{{d}_{{{G}_{1}}{{{\dot{\vee }}}_{T}}{{G}_{2}}}}{{(v)}^{3}}+\sum\limits_{v\in I({{G}_{1}})}{{{d}_{{{G}_{1}}{{{\dot{\vee }}}_{T}}{{G}_{2}}}}{{(v)}^{3}}}}}\\
          &=&\sum\limits_{v\in V({{G}_{1}})}{{{(2{{d}_{{{G}_{1}}}}(v)+{{n}_{2}})}^{3}}}+\sum\limits_{v\in V({{G}_{2}})}{{{({{d}_{{{G}_{2}}}}(v)+{{n}_{1}})}^{3}}}\\
          &&+\sum\limits_{uv\in E({{G}_{1}})}({{{{d}_{{{G}_{1}}}}(u)+{{d}_{{{G}_{1}}}}(v)})^{3}}\\
          &=&\sum\limits_{v\in V({{G}_{1}})}{[8{{d}_{{{G}_{1}}}}{{(v)}^{3}}+12{{n}_{2}}{{d}_{{{G}_{1}}}}(v)^{2}+6{{n}_{2}}^{2}{{d}_{{{G}_{1}}}}(v)+n_{2}^{3}]}\\
          &&+\sum\limits_{v\in V({{G}_{2}})}{[{{d}_{{{G}_{2}}}}{{(v)}^{3}}+3{{n}_{1}}{{d}_{{{G}_{2}}}}(v)^{2}+3{{n}_{1}}^{2}{{d}_{{{G}_{2}}}}(v)+n_{1}^{3}]}\\
          &&+\sum\limits_{uv\in E({{G}_{1}})}[{{d}_{{{G}_{1}}}}(u)^{3}+{{d}_{{{G}_{1}}}}(v)^{3}+3{{d}_{{{G}_{1}}}}(u){{d}_{{{G}_{1}}}}(v)({{d}_{{{G}_{1}}}}(u)+{{d}_{{{G}_{1}}}}(v))]\\
          &=&8\sum\limits_{v\in V({{G}_{1}})}{{{d}_{{{G}_{1}}}}{{(v)}^{3}}}+12{{n}_{2}}\sum\limits_{v\in V({{G}_{1}})}{{{d}_{{{G}_{1}}}}(v)^{2}}+6{{n}_{2}}^{2}\sum\limits_{v\in V({{G}_{1}})}{{{d}_{{{G}_{1}}}}(v)}\\
          &&+{{n}_{1}}n_{2}^{3}+\sum\limits_{v\in V({{G}_{2}})}{{{d}_{{{G}_{2}}}}{{(v)}^{3}}}+3{{n}_{1}}\sum\limits_{v\in V({{G}_{2}})}{{{d}_{{{G}_{2}}}}(v)^{2}}\\
          &&+3{{n}_{1}}^{2}\sum\limits_{v\in V({{G}_{2}})}{{{d}_{{{G}_{2}}}}(v)}+{{n}_{2}}n_{1}^{3}+\sum\limits_{uv\in E({{G}_{1}})}({{d}_{{{G}_{1}}}}(u)^{3}+{{d}_{{{G}_{1}}}}(v)^{3})\\
          &&+3\sum\limits_{uv\in E({{G}_{1}})}{{d}_{{{G}_{1}}}}(u){{d}_{{{G}_{1}}}}(v)({{d}_{{{G}_{1}}}}(u)+{{d}_{{{G}_{1}}}}(v))\\
          &=&8{{F}}({{G}_{1}})+{{F}}({{G}_{2}})+12{{n}_{2}}{{M}_{1}}({{G}_{1}})+3{{n}_{1}}{{M}_{1}}({{G}_{2}})+{{M}_{4}}({{G}_{1}})\\
          &&+3{{ReZM}({{G}_{1}}})+12{{m}_{1}}{{n}_{2}^{2}}+6{{m}_{2}}{{n}_{1}^{2}}+n_{1}{{n}_{2}}({{n}_{1}}^2+n_{2}^{2}).
\end{eqnarray*}
Which is the required result.                                     \qed
\begin{ex}
Applying theorem 4, we get
\begin{eqnarray*}
(i)~{{F}}({{P}_{n}}{{\dot{\vee }}_{T}}{{P}_{m}})&=&mn(m^2+n^2)+6n^2(m-1)+12m^2(n-1)+60mn-64m\\
&&+110n-264,~  m,n\ge 3,\\
(ii)~{{F}}({{P}_{n}}{{\dot{\vee }}_{T}}{{C}_{m}})&=&mn(m^2+n^2)+6n^2m+12m^2(n-1)+60mn-64m\\
&&+128n-250,~   m,n\geq 3,\\
(iii)~{{F}}({{C}_{n}}{{\dot{\vee }}_{T}}{{C}_{m}})&=&mn(m^2+n^2)+6n^2m+12m^2n+60mn+8m\\
&&+128n,~m,n\geq 3,\\
(iv)~{{F}}({{C}_{n}}{{\dot{\vee }}_{T}}{{P}_{m}})&=&mn(m^2+n^2)+6n^2(m-1)+12m^2n+60mn+8m\\
&&+110n-14,~m,n\geq 3.
\end{eqnarray*}
\end{ex}

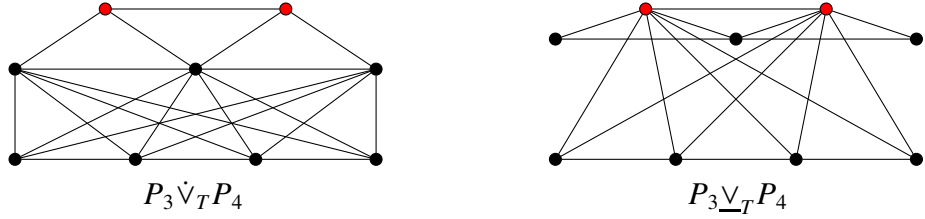
\begin{figure*}[h]
\begin{center}
\begin{tikzpicture}[scale=.4]
\tikzstyle{every node}=[draw, shape=circle, fill=red, scale=.4]
  \node (n1) at (6,7) {};
  \node (n2) at (12,7)  {};
\tikzstyle{every node}=[draw, shape=circle, fill=black, scale=.4]
  \node (n3) at (3,5)  {};
  \node (n4) at (9,5) {};
  \node (n5) at (15,5)  {};
  \node (n6) at (3,2)  {};
  \node (n7) at (7,2)  {};
  \node (n8) at (11,2)  {};
  \node (n9) at (15,2)   {};
  \foreach \from/\to in {n1/n3,n1/n4,n2/n4,n2/n5,n1/n2,n3/n4,n4/n5,n3/n6,n3/n7,n3/n8,n3/n9,n4/n6,n4/n7,n4/n8,n4/n9,n5/n6,n5/n7,n5/n8,n5/n9,n6/n7,n7/n8,n8/n9}
    \draw (\from) -- (\to);
\end{tikzpicture}
\hspace*{2cm}
\begin{tikzpicture}[scale=.4]
\tikzstyle{every node}=[draw, shape=circle, fill=red, scale=.4]
  \node (n1) at (6,7) {};
  \node (n2) at (12,7)  {};
\tikzstyle{every node}=[draw, shape=circle, fill=black, scale=.4]
  \node (n3) at (3,6)  {};
  \node (n4) at (9,6) {};
  \node (n5) at (15,6)  {};
  \node (n6) at (3,2)  {};
  \node (n7) at (7,2)  {};
  \node (n8) at (11,2)  {};
  \node (n9) at (15,2)   {};
  \foreach \from/\to in {n1/n3,n1/n4,n2/n4,n2/n5,n3/n4,n4/n5,n1/n2,n1/n6,n1/n7,n1/n8,n1/n9,n2/n6,n2/n7,n2/n8,n2/n9,n6/n7,n7/n8,n8/n9}
    \draw (\from) -- (\to);
\end{tikzpicture}
{${{P}_{3}}{\dot{\vee}_{T}}{{P}_{4}}$~~~~~~~~~~~~~~~~~~~~~~~~~~~~~~~~~~~~~~~~~~~~~~~~~~~~~~~~${{P}_{3}}{\underline{\vee}_{T}}{{P}_{4}}$}
\caption{The example of vertex T-join and edge T-join of graphs.}
\end{center}
\end{figure*}

\begin{definition}
Let ${{G}_{1}}$ and ${{G}_{2}}$ be two vertex disjoint graphs the edge T-join graph is obtained from $T({{G}_{1}})$ and ${G}_{2}$ by joining each vertex of $T({{G}_{1}})$ to every vertex of ${{G}_{2}}$ and is denoted by ${{G}_{1}}{\underline{\vee}_{T}}{{G}_{2}}$.
\end{definition}

\begin{thm}
If ${{G}_{1}}$ and ${{G}_{2}}$ be two connected graph then
\begin{eqnarray*}
{{F}}({{G}_{1}}{{\underline{\vee }}_{T}}{{G}_{2}})&=&8{{F}}({{G}_{1}})+{{F}}({{G}_{2}})+3{{n}_{2}}^{2}{{M}_{1}}({{G}_{1}})+3{{m}_{1}}{{M}_{1}}({{G}_{2}})+{{M}_{4}}({{G}_{1}})\\
          &&+3{{n}_{2}}{{HM}({{G}_{1}})}+3{{ReZM}({{G}_{1}})}+m_{1}^{2}(6{{m}_{2}}+{{m}_{1}}{{n}_{2}})+{{m}_{1}}{{{n}_{2}}}^{3}.
\end{eqnarray*}
\end{thm}
\n\textit{Proof.} From definition of F-index, we have
\begin{eqnarray*}
{{F}}({{G}_{1}}{{\underline{\vee }}_{T}}{{G}_{2}})&=&\sum\limits_{v\in V({{G}_{1}}{{{\underline{\vee }}}_{T}}{{G}_{2}})}{{{d}_{{{G}_{1}}{{{\underline{\vee }}}_{T}}{{G}_{2}}}}{{(v)}^{3}}}\\
          &=&\sum\limits_{v\in V({{G}_{1}})}{{{d}_{{{G}_{1}}{{{\underline{\vee }}}_{T}}{{G}_{2}}}}{{(v)}^{3}}+\sum\limits_{v\in V({{G}_{2}})}{{{d}_{{{G}_{1}}{{{\underline{\vee }}}_{T}}{{G}_{2}}}}{{(v)}^{3}}+\sum\limits_{v\in I({{G}_{1}})}{{{d}_{{{G}_{1}}{{{\underline{\vee }}}_{T}}{{G}_{2}}}}{{(v)}^{3}}}}}\\
          &=&\sum\limits_{v\in V({{G}_{1}})}{(2{{d}_{{{G}_{1}}}}(v))^{3}}+\sum\limits_{v\in V({{G}_{2}})}{{({{d}_{{{G}_{2}}}}(v)+{{m}_{1}})}^{3}}\\
          &&+\sum\limits_{uv\in E({{G}_{1}})}({{{{d}_{{{G}_{1}}}}(u)+{{d}_{{{G}_{1}}}}(v)+{{n}_{2}})^{3}}}\\
          &=&\sum\limits_{v\in V({{G}_{1}})}{(2{{d}_{{{G}_{1}}}}{{(v)})^{3}}}+\sum\limits_{v\in V({{G}_{2}})}{[{{d}_{{{G}_{2}}}}{{(v)}^{3}}+3{{m}_{1}}{{d}_{{{G}_{2}}}}(v)^{2}}\\
          &&+3{{m}_{1}}^{2}{{d}_{{{G}_{2}}}}(v)+m_{1}^{3}]\\
          &&+\sum\limits_{uv\in E({{G}_{1}})}[({{{{d}_{{{G}_{1}}}}(u)+{{d}_{{{G}_{1}}}}(v))^{3}}}+3{{n}_{2}}{{({{d}_{{{G}_{1}}}}(u)+{{d}_{{{G}_{1}}}}(v))^{2}}}\\
          &&+3{{n}_{2}}^{2}{{({{d}_{{{G}_{1}}}}(u)+{{d}_{{{G}_{1}}}}(v))}}+{{n}_{2}}^{3}]\\
          &=&8\sum\limits_{v\in V({{G}_{1}})}{{{d}_{{{G}_{1}}}}{{(v)}^{3}}}+\sum\limits_{v\in V({{G}_{2}})}{{{d}_{{{G}_{2}}}}{{(v)}^{3}}}+3{{m}_{1}}\sum\limits_{v\in V({{G}_{2}})}{{{d}_{{{G}_{2}}}}(v)^{2}}\\
          &&+3{{m}_{1}}^{2}\sum\limits_{v\in V({{G}_{2}})}{{{d}_{{{G}_{2}}}}(v)}+{{n}_{2}}m_{1}^{3}+\sum\limits_{uv\in E({{G}_{1}})}({{d}_{{{G}_{1}}}}(u)^{3}+{{d}_{{{G}_{1}}}}(v)^{3})\\
          &&+3\sum\limits_{uv\in E({{G}_{1}})}{{d}_{{{G}_{1}}}}(u){{d}_{{{G}_{1}}}}(v)({{d}_{{{G}_{1}}}}(u)+{{d}_{{{G}_{1}}}}(v))+{{n}_{2}}^{3}m_{1}\\
          &&+3{{n}_{2}}\sum\limits_{uv\in E({{G}_{1}})}({{d}_{{{G}_{1}}}}(u)+{{d}_{{{G}_{1}}}}(v))^{2}+3{{n}_{2}}^{2}\sum\limits_{uv\in E({{G}_{1}})}({{d}_{{{G}_{1}}}}(u)+{{d}_{{{G}_{1}}}}(v))\\
          &=&8{{F}}({{G}_{1}})+{{F}}({{G}_{2}})+3{{n}_{2}}^{2}{{M}_{1}}({{G}_{1}})+3{{m}_{1}}{{M}_{1}}({{G}_{2}})+{{M}_{4}}({{G}_{1}})\\
          &&+3{{n}_{2}}{{HM}({{G}_{1}})}+3{{ReZM}({{G}_{1}})}+m_{1}^{2}(6{{m}_{2}}+{{m}_{1}}{{n}_{2}})+{{m}_{1}}{{{n}_{2}}}^{3}.
\end{eqnarray*}
Which is the required result.                                         \qed
\begin{ex}
Using theorem 8, we get
\begin{eqnarray*}
(i)~   {{F}}({{P}_{n}}{{\underline{\vee }}_{T}}{{P}_{m}})&=&m^3(n-1)+3m^2(4n-6)+m(n-1)^3+6(m-1)(n-1)^2\\
&&+60mn-94m+110n-246,~m\geq 2,n\ge 3,\\
(ii) ~  {{F}}({{P}_{n}}{{\underline{\vee }}_{T}}{{C}_{m}})&=&m^3(n-1)+3m^2(4n-6)+m(n-1)^3+6m(n-1)^2\\&&+60mn-94m+128n-250,~m,n\geq 3,\\
(iii) ~ {{F}}({{C}_{n}}{{\underline{\vee }}_{T}}{{C}_{m}})&=&m^3n+12m^2n+mn^2(n+6)+60mn+8m+128n,\\
&&~m,n\geq 3,\\
(iv)  ~ {{F}}({{C}_{n}}{{\underline{\vee }}_{T}}{{P}_{m}})&=&m^3n+12m^2n+mn^3+6n^2(m-1)+60mn\\
&&+8m+110n-14,~m,n\ge 3.
\end{eqnarray*}
\end{ex}

\section{Conclusions}
In this work, we established some useful formula of the F-index of graphs based on the vertex and edge  $\mathcal{F}$-join of graphs where $\mathcal{F}=\{S,R,Q,T\}$, in terms of different topological indices of their factor graphs. Also, we derived some useful examples of the F-index for some important classes of graphs. For further study, still many topological indices for this graph operations can be computed.

\end{document}